\documentclass{amsart}

\usepackage{amssymb,amsfonts}
\usepackage{graphicx}

\usepackage{cite}

\newtheorem{theorem}{Theorem}[section]
\newtheorem{corollary}[theorem]{Corollary}
\newtheorem{proposition}[theorem]{Proposition}

\theoremstyle{definition}
\newtheorem{definition}[theorem]{Definition}

\newtheorem{remark}[theorem]{Remark}

\DeclareMathOperator{\Max}{Max}

\DeclareMathOperator{\Ann}{Ann}

\DeclareMathOperator{\Spec}{Spec}


\numberwithin{equation}{section}

\begin{document}
	
	\title[Some Remarks on Semirings and Their Ideals]{Some Remarks on Semirings and Their Ideals}
	
	\author[Peyman Nasehpour]{\bfseries Peyman Nasehpour}
	
	\address{Peyman Nasehpour\\
		Department of Engineering Science \\
		Golpayegan University of Technology \\
		Golpayegan\\
		Iran}
	\email{nasehpour@gut.ac.ir, nasehpour@gmail.com}
	
	

	\begin{abstract}
		In this paper, we give semiring version of some classical results in commutative algebra related to Euclidean rings, PIDs, UFDs, G-domains, and GCD and integrally closed domains.
	\end{abstract}

	\maketitle
	
	\section{Introduction}
	
	The ideal theoretic method for studying commutative rings has a long and fruitful history \cite{Huckaba1988}. The subject of multiplicative ideal theory is, roughly speaking, the description of multiplicative structure of a ring by means of ideals or certain systems of ideals of that ring \cite{Halter-Koch1998}. Multiplicative ideal theory, originated in the works of Dedekind \cite{Dedekind1894}, Pr\"{u}fer \cite{Prufer1932} and Krull \cite{Krull1968}, is a powerful tool both in commutative algebra \cite{Gilmer1972, LarsenMcCarthy1971} and in a more general context, in commutative monoid theory \cite{Aubert1953, Halter-Koch1998}.
	
	Semirings are ring-like algebraic structures that subtraction is either impossible or disallowed, interesting generalizations of rings and distributive lattices, and have important applications in many different branches of science and engineering \cite{Golan2005}. For general books on semiring theory, one may refer to the resources \cite{Golan1999(a), Golan1999(b), Golan2003, GondranMinoux2008, Glazek2002, HebischWeinert1998}.
	
	Since different authors define semirings differently, it is very important to clarify, from the beginning, what we mean by a semiring in this paper. By a semiring, we understand an algebraic structure, consisting of a nonempty set $S$ with two operations of addition and multiplication such that the following conditions are satisfied:
	
	\begin{enumerate}
		\item $(S,+)$ is a commutative monoid with identity element $0$;
		\item $(S,\cdot)$ is a commutative monoid with identity element $1 \not= 0$;
		\item Multiplication distributes over addition, i.e., $a\cdot (b+c) = a \cdot b + a \cdot c$ for all $a,b,c \in S$;
		\item The element $0$ is the absorbing element of the multiplication, i.e., $s \cdot 0=0$ for all $s\in S$.
	\end{enumerate}
	
	Some of the topics in multiplicative ideal theory for commutative rings have been generalized and investigated for semirings \cite{GhalandarzadehNasehpourRazavi2017, Iseki1956, LaGrassa1995, Lescot2015, Nasehpour2016, Nasehpour2017, NasehpourP2018, Nasehpour2018}. Also see Chapter 7 of the book \cite{Golan1999(b)}.
	
	The purpose of this paper is to generalize some other concepts of multiplicative ideal theory in commutative rings and investigate them in semirings.
	
	Let us recall that a nonempty subset $I$ of a semiring $S$ is said to be an ideal of $S$, if $a+b \in I$ for all $a,b \in I$ and $sa \in I$ for all $s \in S$ and $a \in I$ \cite{Bourne1951}. An ideal $I$ of a semiring $S$ is called a proper ideal of the semiring $S$, if $I \neq S$. A proper ideal $P$ of a semiring $S$ is called a prime ideal of $S$, if $ab\in P$ implies either $a\in P$ or $b\in P$. We collect all the prime ideals of a semiring $S$ in the set $\Spec(S)$. Note that a nonzero and nonunit element $s$ of a semiring $S$ is said to be irreducible if $s = s_1 s_2$ for some $s_1, s_2 \in S$, then either $s_1$ or $s_2$ is unit (multiplicatively invertible). An element $p\in S$ is said to be a prime element, if the principal ideal $(p)$ is a prime ideal of $S$.
	
	Also note that a semiring $S$ is semidomain if $ab = ac$ implies $b=c$ for all $b,c \in S$ and all nonzero $a\in S$. A semidomain $S$ is called to be a unique factorization semidomain (for short UFSD) if the following conditions hold:
	
	\begin{enumerate}
		
		\item[UF1] Each irreducible element of $S$ is a prime element of $S$.
		
		\item[UF2] Any nonzero, nonunit element of $S$ is a product of irreducible elements of $S$.
		
	\end{enumerate}
	
	Section 3 of the paper is devoted to unique factorization semidomains. Beside some basic but interesting results, in Theorem \ref{UFSDKaplansky}, we show that a semidomain $S$ is a UFSD if and only if every nonzero prime ideal of $S$ contains a prime element. The ring version of this result is due to I. Kaplansky \cite{Kaplansky1970}.
	
	Similar to the concept of field of fractions in ring theory, one can define the semifield of fractions $F(S)$ of the semidomain $S$ \cite[p. 22]{Golan1999(a)}.
	
	Note that if $S$ is a semidomain and $F(S)$ its semifield of fractions, then by definition \cite[p. 88]{DulinMosher1972}, an element $u\in F(S)$ is said to be integral over $S$ if there exist $a_1, \ldots, a_n$ and $b_1, \ldots,b_n$ in $S$ such that $u^n+a_1 u^{n-1}+\cdots + a_n = b_1 u^{n-1}+ \cdots + b_n$. The semidomain $S$ is said to be integrally closed if any integral element of $F(S)$ over $S$ belongs to $S$.
	
	We also mention that a nonempty subset $W$ of a semiring $S$ is said to be a multiplicatively closed set (for short an MC-set) if $1\in W$ and for all $w_1,w_2 \in W$, we have $w_1 w_2 \in W$. Note that if $U$ is an MC-set of a semiring $S$, one can define the localization of $S$ at $U$, similar to the definition of the localization in ring theory (refer to \cite{Kim1985} and \cite[\S 11]{Golan1999(b)}). Section 4 is devoted to Integrally closed semirings and similar to ring theory, in Theorem \ref{integrallyclosedsemidomain}, for a semidomain $S$, we prove that the following statements are equivalent:
	
	\begin{enumerate}
		
		\item $S$ is an integrally closed semidomain.
		
		\item $S_T$ is an integrally closed semidomain for any MC-set $T$ of $S$.
		
		\item $S_{\mathfrak{m}}$ is an integrally closed semidomain for any maximal ideal $\mathfrak{m}$ of $S$.
		
	\end{enumerate}
	
	In section 5, we generalize the concept of G-domains \cite[p. 12]{Kaplansky1970} and define a semidomain $S$ to be a Goldman-Krull semidomain if its semifield of fractions $F(S)$ is a finitely generated semiring over $S$. In Corollary \ref{primeintersectionGK}, we show that a semidomain $S$ is a Goldman-Krull semidomain if and only if \[\bigcap_{P\in \Spec(S)-\{0\}} P \neq (0).\]
	
	Let us recall that an ideal $I$ of a semiring $S$ is subtractive if $a+b\in I$ and $a\in I$ imply that $b\in I$ for all $a,b\in S$. An ideal $I$ of a semiring $S$ is called principal if $I = \{sa: s\in S\}$ for some $a\in S$. Such an ideal $I$ of $S$ is denoted by $(a)$. A semiring $S$ is called a principal ideal semidomain (for short PISD) if $S$ is a semidomain and each ideal of $S$ is principal. In Theorem \ref{PISD-GKsemidomain}, we prove that a PISD is a Goldman-Krull semidomain if and only if it has only finitely many prime ideals.
	
	Sections 1 and 2 are devoted to Euclidean and principal ideal semirings. These semirings provide some interesting examples for the other facts that we bring in the rest of the paper.

	\section{Euclidean Semirings}
	
	Let us recall that by definition, a totally ordered set $(O, \leq)$ is said to be well-ordered, if every nonempty subset of $O$ has a least element \cite[p. 38]{Ciesielski1997}. In the following, we clarify what we mean by a Euclidean semiring in this paper:
	
	\begin{definition}
		
		Let $(\mathbb O, \preceq)$ be a well-ordered set with no greatest element. We annex an element $+\infty$ to $\mathbb O$ and denote the new set by $\mathbb O_{\infty}$ and define $x \preceq +\infty$ for any element $x \in \mathbb O_{\infty}$. An $\mathbb O_{\infty}$\emph{-Euclidean norm} on a semiring $S$ is a function $\delta : S \rightarrow \mathbb O_{\infty}$ with the following properties:
		
		\begin{enumerate}
			
			\item $\delta(s) = +\infty$ if and only if $s=0$ for all $s\in S$.
			
			\item If $a,b \in S$ with $b\neq 0$, then there are elements $q,r \in S$ such that $a = bq + r$ with either $r=0$, or $\delta(r) \prec \delta(b)$.
			
		\end{enumerate}
		
		In such a case, we say that $S$ is an $\mathbb O_{\infty}$\emph{-Euclidean semiring}.
		
	\end{definition}
	
	\begin{remark} As far as the author knows, the first resource defining some kinds of Euclidean semirings is the paper \cite{DaleHanson1977}. For more on ``Euclidean semirings'', one can refer to \cite{HebischWeinert1987} and \cite[\S 12]{Golan1999(b)}. We also note that most algebra texts on rings and semirings require an $\mathbb N_{\infty}$-Euclidean norm to have the following additional property:
		
		\begin{itemize}
			\item $\delta (b) \preceq \delta(sb)$, for all $b,s\in S$.
		\end{itemize}
		
		In Proposition \ref{submultiplicative}, which is a generalization of Proposition 12.10 in \cite{Golan1999(b)}, we show that this condition is superfluous.
		
	\end{remark}
	
	\begin{proposition}
		
		\label{submultiplicative}
		
		Let $\delta$ be an $\mathbb O_{\infty}$-Euclidean norm on a semiring $S$. Then, there exists an $\mathbb O_{\infty}$-Euclidean norm $\delta^*$ on $S$ with the following properties:
		
		\begin{enumerate}
			
			\item $\delta^* (a) \preceq \delta (a)$, for all $a\in S$.
			
			\item $\delta^* (b) \preceq \delta(sb)$, for all $b,s\in S$.
			
		\end{enumerate}
		
		\begin{proof}
			Define $\delta^* (a) = \min\{\delta(sa) : s\in S\}$. Clearly, $\delta^*\colon S \rightarrow \mathbb O_{\infty}$ satisfies the conditions (1) and (2). Now, we prove that $\delta^*$ is an $\mathbb O_{\infty}$-Euclidean norm on $S$.
			
			It is straightforward to see that $\delta(s) = +\infty$ if and only if $s=0$ for all $s\in S$. Let $a,b \in S$ with $b\neq 0$. Clearly, there is an $s\in S$, such that $\delta^*(b) = \delta(sb)$. Since $\delta$ is an $\mathbb O_{\infty}$-Euclidean norm on the semiring $S$, there exist elements $q$ and $r$ in $S$ such that $a = qsb+r$, with either $r =0$, or $\delta(r) \prec \delta(sb)$. If $r\neq 0$, we have $\delta^*(r) \preceq \delta(r) \prec \delta(sb)= \delta^*(b)$ and the proof is complete.
		\end{proof}
		
	\end{proposition}

	\begin{proposition}
		
		Let $S$ be an $\mathbb O_{\infty}$-Euclidean semiring. Then, each subtractive ideal of $S$ is principal.
		
		\begin{proof}
			Let $I$ be a nonzero subtractive ideal of $S$ and set $M = \{\delta(s) : s\in I\}$. Obviously $M \subseteq \mathbb O$ and by definition $M$ has a $\preceq$-least element. Let $\delta(b)$ be the $\preceq$-least element of $M$, where $b\in I$. It is clear that the principal ideal $(b)$ is a subset of $I$. Our claim is that $I \subseteq (b)$. Now, take $a\in I-\{0\}$. So, there exist elements $q,r \in S$ such that $a = bq + r$ with either $r=0$ or $\delta(r) \prec \delta(b)$. If $r\neq 0$, then $r\in I$, since $I$ is subtractive and $a,b\in I$ and this means that $\delta(r)$ is smaller than the $\preceq$-least element of $M$, which is obviously a contradiction. Therefore, $r=0$ and $a = bq$ and $a\in (b)$ and this completes the proof.
		\end{proof}
		
	\end{proposition}
	
	A semiring $S$ is said to be a principal ideal semiring (for short PIS) if any ideal of $S$ is principal. A semiring $S$ is called subtractive if each ideal of $S$ is subtractive.
	
	\begin{corollary}
		
		\label{EuclideanPIS}
		
		Each subtractive $\mathbb O_{\infty}$-Euclidean semiring is a PIS.
	\end{corollary}
	
	Let $S$ be a semiring and let $a, b\in S$ such that $b \neq 0$. The element $a$ is said to be a multiple of $b$ if there exists an element $x \in S$ such that $a = bx$. In this case, we say that $b$ divides $a$ or is a divisor of $a$ and write $b \mid a$. Note that this is equivalent to say that $(b) \subseteq (a)$. Also, it is said that $a$ and $b$ are associates if $a=ub$ for some unit $u\in U(S)$ and if $S$ is a semidomain, then this is equivalent to say that $(a) = (b)$.
	
	A \emph{greatest common divisor} of a set $A \subseteq S$, which has at least one nonzero element, is a nonzero element $d$, if $d \mid a$ for any $a\in A$ and if $d^{\prime} \mid a$ for any $a\in A$, then $d^{\prime} \mid d$. It is clear that this is equivalent to say that $(d)$ is the minimal element of all principal ideals containing the ideal generated by the set $A$. A greatest common divisor of the set $A$, which is not necessarily unique, is denoted by $\gcd(A)$. It is easy to see that if the ideal generated by the set $A$ is principal, say the ideal $(a)$, then $a$ is a greatest common divisor of $A$. If at least one of the two $a,b \in S$ is nonzero, then we write $\gcd(a,b)$ for $\gcd(\{a,b\})$.
	
	\begin{proposition}
		
		If $S$ is a semidomain, the set $A \subseteq S$ has at least one nonzero element and $d,d^{\prime}$ are greatest common divisors of $A$, then there exists a unit element $u$ of $S$ such that $d^{\prime} = ud$.
		
		\begin{proof}
			Straightforward.
		\end{proof}
		
	\end{proposition}
	
	\begin{theorem}
		Let $S$ be a subtractive $\mathbb O_{\infty}$-Euclidean semiring and at least one of the two $a,b\in S$ be nonzero. Then, the greatest common divisor $\gcd(a,b)$ of $\{a,b\}$ can be computed algorithmically.
		
		\begin{proof}
			Let for the moment $n\geq 0$ and the two nonzero elements $r_{n-2}, r_{n-1}$ of $S$ are given. Then by definition, there are elements $q_n,r_n \in S$ such that $r_{n-2} = r_{n-1} q_n + r_n$ with either $r_n=0$ or $\delta(r_n) \prec \delta(r_{n-1})$. We can assume that $b \neq 0$ and define $a = r_{-2}$ and $b=r_{-1}$. If each $r_i$ for $i\geq 0$ is nonzero, then we see that $\delta(b) \succ \delta(r_0) \succ \cdots \succ \delta(r_{n-1}) \succ \delta(r_n) \succ \cdots$ is a non-stop descending chain of elements of $\mathbb O$, contradicting this fact that $\mathbb O$ is a well-ordered set. Therefore, there exists an $n\geq0$ such that $r_n = 0$.
			
			Now, let $m$ be the greatest nonnegative number such that $r_{m} \neq 0$. We claim that $r_m = \gcd(a,b)$. Clearly, $r_{m+1} = 0$ and $r_m q_{m+1} = r_{m-1}$. This means that $r_m \mid r_{m-1}$. Also, it is obvious that $r_m \mid r_m$. Using mathematical induction and this fact that $r_{n-2} = r_{n-1} q_n + r_n$ for all $-2 \leq n \leq m$, we get that $r_m \mid b$ and $r_m \mid a$.
			
			On the other hand, $a = bq_0 + r_0$. Since $S$ is subtractive, $r_0$ is an element of the ideal $(a,b)$. But again since $r_{n-2} = r_{n-1} q_n + r_n$ for all $-2 \leq n \leq m$ and $S$ is a subtractive semiring, by mathematical induction, we get that $r_n \in (a,b)$.
		\end{proof}
		
	\end{theorem}
	
	\section{Principal Ideal Semidomains}
	
	We start this section by proving the following proposition for an arbitrary principal ideal semidomain (for short PISD).
	
	\begin{proposition}
		
		\label{primePISD}
		
		Each nonzero prime ideal of a PISD is maximal.
		
		\begin{proof}
			Let $(p)$ be a nonzero prime ideal of the PISD $S$ and let $I=(s)$ be any ideal containing $(p)$. Also, since $p\in (s)$, $p=rs$ for some $r\in S$. Since $(p)$ is prime and $rs\in (p)$, either $r\in (p)$ or $s\in (p)$. If $s\in (p)$, then $(p) = (s) = I$. If not, then there is a $t\in S$ such that $r = pt = rst$. But $S$ is a semidomain and $r\neq 0$, so $st = 1$. This means that $I = S$ and the proof is complete.
		\end{proof}
		
	\end{proposition}
	
	A nonzero, nonunit element $s$ of a semiring $S$ is said to be irreducible if $s = s_1 s_2$ for some $s_1, s_2 \in S$, then either $s_1$ or $s_2$ is a unit. This is equivalent to say that $(s)$ is maximal among proper principal ideals of $S$. An element $p\in S$ is said to be a prime element, if the principal ideal $(p)$ is a prime ideal of $S$, which is equivalent to say if $p \mid ab$, then either $p \mid a$ or $p \mid b$.
	
	\begin{proposition}
		Each prime element of a semidomain is irreducible.
		
		\begin{proof}
			Straightforward.
		\end{proof}
		
	\end{proposition}
	
	\begin{proposition}
		
		\label{irreducibleisprime}
		
		Each nonzero element of a PISD is prime if and only if it is irreducible.
		
		\begin{proof}
			Let $p$ be an irreducible element of $S$. Obviously, the ideal $(p)$ is a proper ideal of $S$ and therefore, it is contained in a maximal ideal $I=(m)$ of $S$. So, there is an element $s\in S$ such that $p = sm$. Since $p$ is irreducible, either $s$ or $m$ is unit. But $m$ cannot be a unit element of $S$, since $I=(m)$ is a maximal ideal of $S$. So, $s$ is a unit and this means that $(p) = (m)$. Finally, any maximal ideal is prime \cite[Corollary 7.13]{Golan1999(b)}. Therefore, $(p)$ is prime and this finishes the proof.
		\end{proof}
		
	\end{proposition}
	
	A semiring $S$ is said to satisfy ACCP property, if any ascending chain of principal ideals of $S$ stops somewhere, i.e., if $$(s_1) \subseteq (s_2) \subseteq \cdots (s_n) \subseteq \cdots \text{~for~} s_1, s_2, \ldots, s_n, \ldots \in S, $$ then there is an $m\in \mathbb N$ such that $(s_m) = (s_{m+i})$ for all $i\geq 1$.
	
	\begin{proposition}
		
		\label{pis-is-accp}
		
		Each PIS satisfies ACCP property.
		
		\begin{proof}
			Take $s_1, s_2, \ldots, s_n, \ldots \in S$ such that $(s_1) \subseteq (s_2) \subseteq \cdots (s_n) \subseteq \cdots $. Set $ I = \bigcup_i (s_i)$. It is easy to see that $I$ is an ideal of $S$ and therefore, since $S$ is a PIS, $I = (a)$ for some $a\in S$. But then, there is an $n\in \mathbb N$ such that $a\in (s_n)$. Now, we have $(a) = (s_n) = (s_{n+i})$ for any $i\geq 1$ and this proves the proposition.
		\end{proof}
	\end{proposition}
	
	Let us recall that if $S$ is a semiring and $f=a_n X^n + \dots + a_1 X + a_0$ is a polynomial in $S[X]$, then the content of $f$, denoted by $c(f)$, is defined to be the finitely generated ideal $(a_n , \dots , a_1, a_0)$ of $S$. Theorem 3 in \cite{Nasehpour2016} states that if $S$ is a subtractive semiring and $f,g\in S[X]$, then there is a natural number $n$ such that $c(f)^n c(g) = c(f)^{n-1} c(fg)$. Also, let us recall that an ideal $I$ of a semiring $S$ is called a cancellation ideal, if $IJ = IK$, implies $J=K$, for all ideals $J,K$ of $S$ \cite[p. 31]{LaGrassa1995}. Clearly, if each nonzero finitely generated ideal of $S$ is a cancellation ideal, then the content formula  $c(f)^n c(g) = c(f)^{n-1} c(fg)$ reduces to $c(fg) = c(f) c(g)$, for all $f,g \in S[X]$. Hence, similar to the definition of Gaussian rings in ring theory, it is natural to define Gaussian semirings as follows:
	
	\begin{definition}
		
		\label{Gaussiansemiring}
		
		A semiring $S$ is said to be Gaussian if $c(fg) = c(f)c(g)$, for all $f,g \in S[X]$ (\cite[Definition 7]{Nasehpour2016}).
		
	\end{definition}
	
	\begin{proposition}
		
		\label{PISD-Gaussian}
		
		If $S$ is a subtractive PISD, then $S$ is Gaussian.
		
		\begin{proof}
			By Theorem 3 in \cite{Nasehpour2016}, if $f,g \in S[X]$, then there is a natural number $n$ such that $$c(f)^n c(g) = c(f)^{n-1} c(fg).$$ On the other hand, each ideal of the semidomain $S$ is principal. So, each nonzero ideal of $S$ is a cancellation ideal of $S$ and therefore, $c(fg) = c(f)c(g)$, for all $f,g \in S[X]$ and this means that the semiring $S$ is Gaussian and the proof is complete.
		\end{proof}
	\end{proposition}
	
	Clearly, by Corollary \ref{EuclideanPIS}, we have the following result:
	
	\begin{corollary}
		
		\label{EuclideanGaussian}
		
		Each subtractive $\mathbb O_{\infty}$-Euclidean semidomain is Gaussian.
		
	\end{corollary}

	\section{Unique Factorization Semidomains}
	
	A semiring $S$ is called a unique factorization semidomain (for short UFSD) or sometimes factorial semiring, if the following conditions are satisfied:
	
	\begin{enumerate}
		
		\item[UF1] Each irreducible element of $S$ is a prime element of $S$.
		
		\item[UF2] Any nonzero, nonunit element of $S$ is a product of irreducible elements of $S$.
		
	\end{enumerate}
	
	\begin{proposition}
		
		\label{accp1}
		
		Let $S$ be a semidomain such that it satisfies ACCP property. Any nonzero and nonunit element of $S$ factors into a product of irreducible elements of $S$.
		
		\begin{proof}
			Let $s$ be a nonzero and nonunit element of the semiring $S$ such that it cannot be factored into irreducible elements of $S$. Therefore, $s$ is not irreducible and there are $s_1, t_1 \in S$ such that $s = s_1 t_1$ and $(s) \subset (s_1)$ and $(s) \subset (t_1)$. But one of the two elements $s_1$ and $t_1$ cannot be factored into irreducible elements of $S$, since $s$ cannot be done the same. Now let, for example, $s_1$ does not have a factorization into irreducibles. By iterating the same argument, we may construct an ascending chain of principals $(s) \subset (s_1) \subset \cdots \subset (s_n) \subset \cdots $, which does not stop anywhere. But this contradicts the hypothesis of $S$ having the property ACCP and finishes the proof.
		\end{proof}
		
\end{proposition}
	
\begin{theorem}
		
		\label{PISDisUFSD}
		
		Each PISD is UFSD.
		
\begin{proof}
Let the semiring $S$ be a PISD. By Proposition \ref{irreducibleisprime}, any irreducible element of $S$ is a prime element of $S$. Also, by Proposition \ref{pis-is-accp}, $S$ has ACCP property. Now, let $s$ be a nonzero, nonunit element of $S$. So, by Proposition \ref{accp1}, $s$ is factored into irreducible elements of $S$. Therefore, we have already proved that $S$ is UFSD.
\end{proof}
		
\end{theorem}
	
	\begin{corollary}
		Each subtractive $\mathbb O_{\infty}$-Euclidean semidomain is a UFSD.
	\end{corollary}
	
	A set-theoretic complement $W=S-P$ of a prime ideal $P$ has this property that $ab \in W$ if and only if $a,b\in W$ for all $a,b\in S$. We give a name to this property in the following definition:
	
	\begin{definition}
		
		Let $S$ be a semiring. An MC-set $W$ of a semiring $S$ is called \emph{saturated} if $ab \in W$ if and only if $a,b\in W$ for all $a,b\in S$.
		
	\end{definition}
	
	\begin{theorem}
		
		\label{saturatedunionofprimes}
		
		Let $S$ be a semiring and $W$ a nonempty subset of $S$. Then, the following statements are equivalent:
		
		\begin{enumerate}
			
			\item $W$ is a saturated MC-set of $S$.
			
			\item The complement of $W$ is a union of prime ideals of $S$.
			
		\end{enumerate}
		
		\begin{proof}
			$(2) \Rightarrow (1)$: Straightforward.
			
			$(1) \Rightarrow (2)$: Let $s$ be in the complement of $W$. Since $W$ is saturated, $(s)$ is disjoint from $W$. By Theorem 2.1 in \cite{Nasehpour2017}, $(s)$ is a subset of a prime ideal of $S$, which is also disjoint from $W$. This means that any element of $W$ is an element of a prime ideal of $S$ such that it is disjoint from $W$ and this completes the proof.
		\end{proof}
		
	\end{theorem}
	
	\begin{proposition}
		
		\label{saturatedKaplansky}
		
		Let $S$ be a semidomain and $W$ the set of all elements of $S$ expressible as a finite product of a unit and some prime elements of $S$. Then $W$ is a saturated MC-set.
		
		\begin{proof}
			Since empty product is supposed to be $1$, $1\in W$. Now it is clear that $W$ is an MC-set.  In the next step we prove that $W$ is saturated. Our proof is by induction on the number of prime factors of the elements of $W$. Take $ab\in W$. If $ab$ has no factor of  prime elements, then it needs to be unit. So, $a$ and $b$ are both unit, which implies that $a,b \in W$. If $ab = p$, where $p$ is a prime element of $W$, then we may suppose that $p$ divides $a$, which means that $a=p a_1$. This implies that $b$ is a unit and so $a,b \in W$. Finally, let $ab = p_1 \cdots p_n$. Obviously, $p_1$ divides one of the two $a$ and $b$, say $a$. So, $a=p_1 a_1$ and therefore, $a_1 b = p_2 \cdots p_n$. By induction, both $a_1$ and $b$ are in $W$, which implies that $a$ is also in $W$ and the proof is complete.
		\end{proof}
		
	\end{proposition}
	
	The ring version of the following theorem is due to the Canadian mathematician Irving Kaplansky (1917--2006).
	
	\begin{theorem}
		
		\label{UFSDKaplansky}
		
		A semidomain $S$ is a UFSD if and only if every nonzero prime ideal of $S$ contains a prime element.
		
		\begin{proof}
			$(\Rightarrow)$: Let $S$ be a UFSD and $P$ a nonzero prime ideal of $S$. Take $a\in P$ to be nonzero and nonunit. Then obviously, $a=u p_1 \cdots p_n$ such that $u$ is a unit and $p_i$ is a prime element of $S$ for any $i$. Since $P$ is a prime ideal of $S$, there is an $i$ such that $p_i \in P$.
			
			$(\Leftarrow)$: Suppose that $S$ is a semidomain and each nonzero prime ideal of $S$ contains a prime element. Take $W$ to be the set of all elements of $S$ expressible as a finite product of a unit and some prime elements of $S$. Imagine $s$ is nonzero and not an element of $W$. Since by Proposition \ref{saturatedKaplansky}, $W$ is a saturated MC-set, the principal ideal $(s)$ is disjoint from $W$. By Theorem 2.1 in \cite{Nasehpour2017}, $(s)$ can be expanded to a prime ideal $P$ of $S$, which is disjoint from $W$. By hypothesis, $P$ contains a prime element, a contradiction. Q.E.D.
		\end{proof}
		
	\end{theorem}
	
	\section{Integrally Closed Semidomains}
	
	Similar to the concept of field of fractions in ring theory, one can define the semifield of fractions $F(S)$ of the semidomain $S$ \cite[p. 22]{Golan1999(a)}.
	
	\begin{definition}
		
		\label{defintegrallyclosed}
		
		(\cite[p. 88]{DulinMosher1972}) Let $S$ be a semidomain and $F(S)$ its semifield of fractions. The element $u\in F(S)$ is said to be integral over $S$ if there exist $a_1, \ldots, a_n$ and $b_1, \ldots,b_n$ in $S$ such that $u^n+a_1 u^{n-1}+\cdots + a_n = b_1 u^{n-1}+ \cdots + b_n$. The semidomain $S$ is said to be integrally closed if any integral element of $F(S)$ over $S$ belongs to $S$.
		
	\end{definition}
	
	Let us recall that if $U$ is an MC-set of a semiring $S$, one can define the localization of $S$ at $U$, similar to the definition of the localization in ring theory. For a good introduction to the concept of localization of semirings, one can refer to \cite{Kim1985} and also \cite[\S 11]{Golan1999(b)}.
	
	\begin{theorem}
		
		\label{integrallyclosedsemidomain}
		
		Let $S$ be a semidomain. Then the following statements are equivalent:
		
		\begin{enumerate}
			
			\item $S$ is an integrally closed semidomain.
			
			\item $S_T$ is an integrally closed semidomain for any MC-set $T$ of $S$.
			
			\item $S_{\mathfrak{m}}$ is an integrally closed semidomain for any maximal ideal $\mathfrak{m}$ of $S$.
			
		\end{enumerate}
		
		\begin{proof}
			
			Let $S$ be a semidomain and $F(S)$ its semifield of fractions.
			
			$(1) \Rightarrow (2)$: Let $T$ be an MC-set in $S$. Obviously, one can consider $S_T$ as a subsemiring of $F(S)$. Therefore, $S_T$ is a semidomain. Now, take $u$ to be an element in the semifield of fractions of $S_T$. This implies that $u$ satisfies in the following equation:
			$$(A): u^n+a_1/t_1 u^{n-1}+\cdots + a_n / t_n = b_1 /t^{\prime}_1 u^{n-1}+ \cdots + b_n /t^{\prime}_n ,$$
			
			where $a_i,b_i \in S$ and $t_i, t^{\prime}_i \in T$.\\
			
			Define $t= (t_1 \cdots t_n)(t^{\prime}_1 \cdots t^{\prime}_n)$, $v_i = t/t_i$ and $v^{\prime}_i = t / t^{\prime}_i$. Now, if we multiply the equation $(A)$ by $t$, we have the following equation: $$ (B): tu^n+a_1 v_1 u^{n-1}+\cdots + a_n v_n = b_1 v^{\prime}_1 u^{n-1}+ \cdots + b_n v^{\prime}_n .$$
			
			Finally, if we multiply the equation $(B)$ by $t^{n-1}$, we get an equation, which says that $tu$ is integral over $S$. Since $S$ is integrally closed, $tu \in S$, which implies that $u\in S_T$.
			
			$(2) \Rightarrow (3)$: Trivial.
			
			$(3) \Rightarrow (1)$: Let $S_{\mathfrak{m}}$ be an integrally closed semidomain for any maximal ideal $\mathfrak{m}$ of $S$. The semidomain $\bigcap_{\mathfrak{m} \in \Max(S)} S_{\mathfrak{m}}$ is also an integrally closed semidomain. Our claim is that $\bigcap_{\mathfrak{m} \in \Max(S)} S_{\mathfrak{m}} = S$. It is obvious that $S \subseteq S_{\mathfrak{m}}$ for all $\mathfrak{m} \in \Max(S)$.
			
			Now, we prove that $\bigcap_{\mathfrak{m} \in \Max(S)} S_{\mathfrak{m}} \subseteq S$. Take $y\in F(S)-S$. Define $I := \{x\in S : xy \in S\}$. It is easy to check that $I$ is an ideal of $S$ and since $y\notin S$, $I\neq S$. Therefore, $I$ is a subset of a maximal ideal, say $\mathfrak{m}$, of $S$. Our claim is that $y\notin S_{\mathfrak{m}}$. On the contrary, assume that $y\in S_{\mathfrak{m}}.$ This implies that $sy\in S$ for some $s\in S-\mathfrak{m}$. So, $s\in I$, a contradiction and this finishes the proof.
		\end{proof}
		
	\end{theorem}
	
	In the proof of Theorem \ref{integrallyclosedsemidomain}, we used the equality $\bigcap_{\mathfrak{m} \in \Max(S)} S_{\mathfrak{m}} = S$. Now, we apply a similar technique to prove that if a semiring is locally nilpotent-free, then it is itself nilpotent-free. We recall that a semiring $S$ is nilpotent-free if the assumption $s\in S-\{0\}$ implies that $s^n \neq 0$, for all natural numbers $n$.
	
	\begin{theorem}
		
		Let $S$ be a semiring such that for any maximal ideal $\mathfrak{m}$, $S_\mathfrak{m}$ is nilpotent-free. Then $S$ is nilpotent-free.
		
		\begin{proof}
			One can easily check that in order to prove that a semiring $S$ is nilpotent-free, it is enough to prove that if $s^2 = 0$, then $s = 0$ for any $s\in S$.
			
			Now, let $s\in S-\{0\}$ with $s^2 =0$. It is clear that $\Ann(s)$ is a proper ideal of $S$ and is, therefore, contained in a maximal ideal $\mathfrak{m}$ of $S$. Our claim is that $s/1$ is nonzero in $S_\mathfrak{m}$. For, if $s/1$ is zero in $S_\mathfrak{m}$, then $ts =0$ for some $t\notin \mathfrak{m}$ and this implies that $t$ is an annihilator of $s$ and therefore is an element of $\mathfrak{m}$, a contradiction. Now, it is obvious that $(s/1)^2$ is zero in $S_\mathfrak{m}$, which means that $S_\mathfrak{m}$ is not nilpotent-free and the proof is complete.
		\end{proof}
		
	\end{theorem}
	
	\begin{definition}
		
		A semidomain $S$ is said to be a GCD semidomain if $\gcd(a,b)$ exists for any $a,b\in S$, whenever at least one of the elements $a$ and $b$ is nonzero.
		
	\end{definition}
	
	\begin{remark}
		
		\begin{enumerate}
			
			\item Any unique factorization semidomain is a GCD semidomain.
			
			\item A semiring is said to be a B\'{e}zout semiring if any 2-generated ideal of $S$ is principal. It is clear that any B\'{e}zout semidomain is a GCD semidomain.
			
		\end{enumerate}
		
	\end{remark}
	
	\begin{proposition}
		
		\label{gcdequalities}
		
		Let $S$ be a GCD semidomain. Then the following statements hold for all $a,b,c,d\in S$:
		
		\begin{enumerate}
			
			\item $\gcd(ab,ac) = a \gcd(b,c)$.
			
			\item If $\gcd(a,b) = d$ then $\gcd(a/d,b/d) = 1$.
			
			\item If $\gcd(a,b) = 1$ and $\gcd(a,c) = 1$ then $\gcd(a,bc) = 1$.
			
		\end{enumerate}
		
		\begin{proof}
			$(1)$: Let $h = \gcd(ab,ac)$. Since $a$ divides both $ab$ and $ac$, $a$ divides $h$. So, $h = a g$. But $ag$ divides both $ab$ and $ac$, which means that $g$ divides both $b$ and $c$. If $f$ divides $b$ and $c$, then $af$ divides $ab$ and $ac$, which means that $af$ divides $h = ag$. But this implies that $f$ divides $g$. This proves that $g = \gcd(b,c)$.
			
			$(2)$ is an immediate consequence of $(1)$.
			
			$(3)$: Let $d$ divides $a$ and $bc$. So, $d$ divides both $ab$ and $bc$, which implies that $d$ divides $\gcd(ab,bc)$. But by $(1)$, $\gcd(ab,bc) = b \gcd(a,c) = b $. So $d$ divides both $a$ and $b$, which means that $\gcd(a,bc) = 1$.
		\end{proof}
		
	\end{proposition}
	
	Now we prove that many GCD semidomains are integrally closed:
	
	\begin{theorem}
		
		\label{gcdisintegrallyclosed}
		
		Let $S$ be a GCD semidomain such that each principal ideal of $S$ is subtractive. Then, $S$ is integrally closed.
		
		\begin{proof}
			
			Let $S$ be a GCD semidomain and $F(S)$ its semifield of fractions. Suppose $u\in F(S)$ such that it satisfies the following equation:
			
			$$ u^n+a_1 u^{n-1}+\cdots + a_n = b_1 u^{n-1}+ \cdots + b_n. $$
			
			Imagine $u = s/t$. We can suppose from the beginning that $\gcd(s,t) =1$, since if $\gcd(s,t) = d$, then $u = (s/d) / (t/d)$, while $\gcd(s/d,t/d) =1$. Now, we have the following equation:
			$$ s^n+a_1 s^{n-1} t+\cdots + a_n t^n = b_1 s^{n-1} t+ \cdots + b_n t^n. $$
			
			But by assumption, any principal ideal of $S$ is subtractive. So, $s^n \in (t)$, which means that $t$ divides $s^n$. But by Proposition \ref{gcdequalities}, $\gcd(t,s^n) = 1$, which implies that $t$ is a unit and $u\in S$, as required.
		\end{proof}
		
	\end{theorem}
	
	\section{Goldman-Krull Semidomains}
	
	Let us recall that if $T$ is a semiring and $S$ is a subset of $T$, then $S$ is said to be a subsemiring of the semiring $T$, if $a+b, ab \in S$ for all $a,b \in S$ and $0_T , 1_T \in S$.
	
	Let $S$ be a subsemiring of the semiring $T$. $T$ is said to be a finitely generated semiring over $S$, if there are a finite number of the elements $t_1, \ldots, t_n$ of $T$ such that any element $x$ of $T$ can be of the form $x= s_1 x_1 + \cdots + s_m x_m$, where $s_1, \ldots,s_m$ belong to $S$ and $x_1 , \ldots , x_m$ belong to the multiplicative submonoid of $T$ generated by $t_1, \ldots, t_n$. In such a case, we write that $T = S[t_1, \ldots,t_n]$.
	
	\begin{proposition}
		
		\label{GoldmanKrull}
		Let $S$ be a semidomain and $F(S)$ its semifield of fractions. Then the following statements are equivalent:
		
		\begin{enumerate}
			
			\item $F(S)$ is a finitely generated semiring over $S$.
			
			\item $F(S)$ can be generated by a single element of $S$.
			
		\end{enumerate}
		
		\begin{proof}
			$(2) \Rightarrow (1)$ is obvious. For $(1) \Rightarrow (2)$, take $F(S) = S[s_1/u_1, \ldots, s_n/u_n]$. It is clear that if we set $u=u_1 \cdots u_n$, we have $F(S)=S[1/u]$.
		\end{proof}
		
	\end{proposition}
	
	\begin{remark}
		Domains satisfying the equivalent conditions in Proposition \ref{GoldmanKrull}, are called G-domains \cite[p. 12]{Kaplansky1970}. G-domains appeared in the works of American mathematician Oscar Goldman (1925--1986) and German mathematician Wolfgang Krull (1899--1971) (cf. \cite{Goldman1951} and \cite{Krull1951}). This justifies us to give the following definition:
	\end{remark}
	
	\begin{definition}
		
		\label{GKsemidomaindef}
		
		We define a semidomain $S$ to be a Goldman-Krull semidomain, if its semifield of fractions $F(S)$ is a finitely generated semiring over $S$.
	\end{definition}
	
	\begin{theorem}
		
		\label{GKsemidomain}
		
		Let $S$ be a semidomain and $F(S)$ its semifield of fractions. Then the following statements for a nonzero element $u$ of $S$ are equivalent:
		
		\begin{enumerate}
			
			\item Any nonzero prime ideal of $S$ contains $u$.
			\item Any nonzero ideal of $S$ contains a power of $u$.
			\item $F(S) = S[1/u]$.
			
		\end{enumerate}
		
		\begin{proof}
			
			$(1) \Rightarrow (2)$: Let $I$ be a nonzero ideal of $S$ such that it doesn't contain any power of $u$. By Theorem 2.1 in \cite{Nasehpour2017}, $I$ can be expanded to a prime ideal $P$ disjoint from the MC-set $\{u^n\}^{\infty}_{n=1}$, contradicting our assumption.
			
			$(2) \Rightarrow (3)$: Let $s$ be a nonzero element of $S$. The principal ideal $(s)$ contains a power of $u$, which means that there is a nonzero $t\in S$ such that $u^n = st$. This implies that $s^{-1} = t u^{-n} \in S[1/u]$. It is now clear that $F(S)=S[1/u]$.
			
			$(3) \Rightarrow (1)$: Let $P$ be a nonzero prime ideal of $S$ and take $s\in P$ to be a nonzero element. Since $F(S)=S[1/u]$, $s^{-1} = t u^{-n}$ for some $t\in S$ and $n\in \mathbb N$. This implies that $u^n = st \in P$ and finally $u\in P$ and the proof is complete.
		\end{proof}
		
	\end{theorem}
	
	\begin{corollary}
		Let $S$ be a Goldman-Krull semidomain and $T$ be a semidomain such that $S \subseteq T \subseteq F(S)$. Then, $T$ is a Goldman-Krull semidomain.
		
		\begin{proof}
			Let $S$ and $T$ be semidomains such that $S \subseteq T \subseteq F(S)$. Clearly, the semifield of fractions $F(S)$ of $S$ is the same as the semifield of fractions $F(T)$ of $T$. But $S$ is a  Goldman-Krull semidomain. So, by Theorem \ref{GKsemidomain}, there is a nonzero $u\in S$ such that $F(S)= S[1/u]$. This implies that $F(T) = T[1/u]$ and the proof is complete.
		\end{proof}
		
	\end{corollary}
	
	\begin{corollary}
		
		\label{primeintersectionGK}
		
		A semidomain $S$ is a Goldman-Krull semidomain if and only if \[\bigcap_{P\in \Spec(S)-\{0\}} P \neq (0).\]
		
		\begin{proof}
			By Theorem \ref{GKsemidomain}, $S$ is a Goldman-Krull semidomain if and only if it has a nonzero element $u$ lying in every nonzero prime ideal of $S$.
		\end{proof}
		
	\end{corollary}
	
	\begin{theorem}
		
		\label{PISD-GKsemidomain}
		
		A PISD is a Goldman-Krull semidomain if and only if it has only finitely many prime ideals.
		
		\begin{proof}
			Let $S$ be a Goldman-Krull PISD. It is clear that by Corollary \ref{primeintersectionGK}, \[I = \bigcap_{P\in \Spec(S)-\{0\}} P \neq (0).\] Suppose that $S$ has infinitely many prime ideals, say, $(p_1), (p_2), \ldots$. If $a\in I$, then $p_i | a$, for each $i$. But by Theorem \ref{PISDisUFSD}, each PISD is a UFSD and therefore, $a = p^{e_1}_1 \cdots p^{e_n}_n$, a contradiction.
			
			On the other hand, if $S$ has finitely many nonzero prime ideals $(p_1), (p_2), \ldots, (p_n)$, then $p_1 p_2 \cdots p_n$ is a nonzero element of $S$ in the intersection of all nonzero prime ideals of $S$. Therefore, $S$ is a Goldman-Krull semidomain by Corollary \ref{primeintersectionGK} and the proof is complete.
		\end{proof}
		
	\end{theorem}
	
	\section*{Acknowledgements}
	
	The author is supported in part by the Department of Engineering Science at the Golpayegan University of Technology and wishes to thank the department for supplying all necessary facilities in pursuing this research. The author is also grateful to Professor Dara Moazzami for his help and encouragements.
	
	\bibliographystyle{plain}

\begin{thebibliography}{15.}
		
		\bibitem{Aubert1953} K. E. Aubert, {\em On the ideal theory of commutative semi-groups}, Math. Scad., \textbf{1} (1953), 39--54.
		
		\bibitem{Bourne1951} S. Bourne, {\em The Jacobson radical of a semiring}, Proc. Nat. Acad. Sci. {\bf 37} (1951), 163--170.
		
		\bibitem{Ciesielski1997} K. Ciesielski, {\em Set Theory for the Working Mathematician}, London Mathematical Society Student Texts {\bf 39}, Cambridge University Press, Cambridge, 1997.
		
		\bibitem{DaleHanson1977} L. Dale and D. L. Hanson, {\em The structure of ideals in a Euclidean semiring}, Kyungpook Math. J. {\bf 17} (1977), 21--29.
		
		\bibitem{Dedekind1894} R. Dedekind, {\em \"{U}ber die Theorie der ganzen algebraiscen Zahlen}, Supplement XI to P.G. Lejeune Dirichlet: Vorlesung \"{u}ber Zahlentheorie 4 Aufl., Druck und Verlag, Braunschweig, 1894.
		
		\bibitem{DulinMosher1972} D. J. Dulin and J. R. Mosher, {\em The Dedekind property for semirings}, J. of Aust. Math. Soc., {\bf 14} (1972), 82--90.
		
		\bibitem{FontanaHuckabaPapick1997} M. Fontana, J. A. Huckaba, and I. J. Papick, {\em Pr\"{u}fer Domains}, Marcel Dekker, New York, 1997.
		
		\bibitem{GhalandarzadehNasehpourRazavi2017} S. Ghalandarzadeh, P. Nasehpour, R. Razavi, {\em Invertible ideals and Gaussian semirings}, Arch. Math. (Brno), {\bf 53} (2017), 179--192.
		
		\bibitem{Gilmer1972} R. Gilmer, {\em Multiplicative Ideal Theory}, Marcel Dekker, New York, 1972.
		
		\bibitem{Glazek2002} K. G{\l}azek, {\em A guide to the literature on semirings and their applications in mathematics and information sciences}, Kluwer, Dordrecht, 2002.
		
		\bibitem{Golan1999(a)} J. S. Golan, {\em Power Algebras over Semirings, With Applications in Mathematics and Computer Science}, Kluwer, Dordrecht, 1999.
		
		\bibitem{Golan1999(b)} J. S. Golan, {\em Semirings and Their Applications}, Kluwer, Dordrecht, 1999.
		
		\bibitem{Golan2003} J. S. Golan, {\em Semirings and Affine Equations over Them: Theory and Applications}, Kluwer, Dordrecht, 2003.
		
		\bibitem{Golan2005} J. S. Golan, {\em Some recent applications of semiring theory}, International Conference on Algebra in Memory of Kostia Beider at National Cheng Kung University, Tainan, 2005.
		
		\bibitem{Goldman1951} O. Goldman, {\em Hilbert rings and the Hilbert Nullstellensatz}, Math. Z. {\bf54} (1951), (2), 136--140.
		
		\bibitem{GondranMinoux2008} M. Gondran, and M. Minoux, {\em Graphs, Dioids and Semirings}, Springer, New York, 2008.
		
		\bibitem{Halter-Koch1998} F. Halter-Koch, {\em Ideal Systems. An Introduction to Multiplicative Ideal Theory}, Marcel Dekker, New York, 1998.
		
		\bibitem{HebischWeinert1987} U. Hebisch and H. J. Weinert, {\em On Euclidean semirings}, Kyungpook Math. J., {\bf 27} (1987), 61--88.
		
		\bibitem{HebischWeinert1998} U. Hebisch and H. J. Weinert, {\em Semirings, Algebraic Theory and Applications in Computer Science}, World Scientific, Singapore, 1998.
		
		\bibitem{Huckaba1988} J. A. Huckaba, {\em Commutative Rings with Zero Divisors}, Marcel Dekker, New York, 1988.
		
		\bibitem{Iseki1956} K. Is\'{e}ki, {\em Ideal theory of semiring}, Proc. Japan. Acad., {\bf 32} (1956), 554-559.
		
		\bibitem{Kaplansky1970} I. Kaplansky, {\em Commutative Rings}, Allyn and Bacon, Boston, 1970.
		
		\bibitem{Kim1985} C. B. Kim, {\em A note on the localization in semirings}, Journal of Scientific Institute at Kookmin Univ., {\bf 3} (1985), 13--19.
		
		\bibitem{Krull1968} W. Krull, {\em Idealtheorie, Zweite, erg\"{a}nzte Auflage}, Springer-Verlag, Berlin, 1968.
		
		\bibitem{Krull1951} W. Krull, {\em Jacobsonsche Ringe, Hilbertscher Nullstellensatz, Dimensionstheorie}, Math. Z. {\bf 54} (4) (1951), 354--387.
		
		\bibitem{LaGrassa1995} S. LaGrassa, {\em Semirings: Ideals and Polynomials}, PhD Thesis, University of Iowa, 1995.
		
		\bibitem{LarsenMcCarthy1971} M. D. Larsen and P. J. McCarthy, {\em Multiplicative Theory of Ideals}, Academic Press, New York, 1971.
		
		\bibitem{Lescot2015} P. Lescot, {\em Prime and primary ideals in semirings}, Osaka J. Math. {\bf 52} (2015), 721--736.
		
		\bibitem{Nasehpour2016} P. Nasehpour, {\em On the content of polynomials over semirings and its applications}, J. Algebra Appl., {\bf 15}, No. 5 (2016), 1650088 (32 pages).
		
		\bibitem{NasehpourP2018} P. Nasehpour, {\em On zero-divisors of semimodules and semialgebras}, arXiv:1702.00810 (2018).
		
		\bibitem{Nasehpour2017} P. Nasehpour, {\em Pseudocomplementation and minimal prime ideals in semirings}, arXiv:1703.08923 (2017), to appear in Algebra Univers.
		
		\bibitem{Nasehpour2018} P. Nasehpour, {\it Valuation semirings},  J. Algebra. Appl, {\bf 16}, No. 4 (2018) 1850073 (23 pages).
		
		\bibitem{Prufer1932} H. Pr\"{u}fer, {\em Untersuchungen \"{u}ber Teilbarkeitseigenschaften in K\"{o}rpern}, J. Reine Angew. Math. {\bf 168} (1932), 1--36.
		
	\end{thebibliography}

\end{document}